\documentclass[11pt,oneside]{amsart}

\title[The conjugacy problem for $\aut\Gamma$]{The conjugacy problem
  for the automorphism group of the random graph}
\author{Samuel Coskey}
\address{Samuel Coskey, Department of Mathematics, The Graduate Center
  of The City University of New York, 365 Fifth Avenue, New York, NY
  10016}
\email{scoskey@nylogic.org}
\author{Paul Ellis}
\address{Paul Ellis, Department of Mathematics, University of
  Connecticut, 196 Auditorium Road Unit 3009, Storrs, CT 06269}
\email{ellis@math.uconn.edu}
\author{Scott Schneider}
\address{Scott Schneider, Mathematics and Computer Science Department,
  Wesleyan University, Science Tower 655, 265 Church St., Middletown,
  CT 06459}
\email{sschneider01@wesleyan.edu}
\thanks{This work was partially supported by NSF grant DMS 0600940.}

\usepackage[marginratio=1:1]{geometry}
\usepackage{setspace}
\usepackage{mathpazo,bm}
\usepackage{mydefs}

\DeclareMathOperator{\Mod}{Mod}

\begin{document}
\begin{abstract}
  We prove that the conjugacy problem for the automorphism group of
  the random graph is Borel complete, and discuss the analogous
  problem for some other countably categorical structures.
\end{abstract}

\maketitle

\section{Introduction}

In \cite{T1}, Truss proved that the automorphism group of the random
graph is simple.  Truss continued to study the conjugacy relation on
this group in \cite{T2}, proving that any element can be written as
the product of at most three conjugates of any other.  He did not,
however, give a complete solution to the conjugacy problem for this
group.  In this paper we shall show in a precise sense that this
problem is in fact as complex as it can conceivably be.  We would
like to thank Simon Thomas for suggesting this project.

In order to make the notion of the complexity of a classification
precise, we will use the theory of definable equivalence relations.  A
\emph{standard Borel space} is a Polish space equipped just with its
$\sigma$-algebra of Borel sets.  For instance, $\RR^n$, Cantor space
$2^\NN$, and the Baire space $\NN^\NN$ are all standard Borel spaces.
Moreover, it is well-known that any Borel subset of a standard Borel
space is again a standard Borel space in its own right.  If $X$ is a
standard Borel space, an equivalence relation $E$ on $X$ is called
\emph{Borel} (analytic, \emph{etc.})\ if it is Borel (analytic,
\emph{etc.})\ as a subset of $X\times X$.  It turns out that
classification problems from many areas of mathematics may be realized
as definable equivalence relations on suitably chosen standard Borel
spaces.

For example, consider the problem of classifying all countable graphs
up to graph isomorphism.  Letting $\mathcal{G}$ be the set of binary
relations $R$ on $\NN$ such that the structure $(\NN,R)$ is an
undirected graph, it is easily checked that $\mathcal{G}$ is a Borel
subset of $2^{\NN\times\NN}$ and hence is itself a standard Borel
space.  More generally, suppose that $\mathcal{L}$ is a countable
relational language and that $\sigma$ is a sentence of the infinitary
logic $\mathcal{L}_{\omega_1,\omega}$ in which countable conjunctions
and disjunctions are allowed.  Then
\[\Mod(\sigma)=\{\mathcal{M}\mid
\textrm{the domain of $\mathcal{M}$ is $\NN$, and }\mathcal{M}\models\sigma\}
\]
is a standard Borel space, and the isomorphism relation $\oiso_\sigma$
on $\Mod(\sigma)$ is an analytic equivalence relation (see for
instance \cite{HK}).

If $E$ is an equivalence relation on the standard Borel space $X$, we
say that $E$ is \emph{smooth}, or concretely classifiable, if there
exists a Borel function $f$ from $X$ into some standard Borel space
$Y$ such that for all $x,x'\in X$,
\[x\mathrel{E}x'\iff f(x)=f(x')\;.
\]
In other words, $f$ selects elements of $Y$ as complete invariants for
the classification of elements of $X$ up to $E$.  For instance, let
$\mathcal C_\NN$ denote the conjugacy relation on $S_\infty$
(\emph{i.e.}, the set $\aut\NN$ of permutations $f\from\NN\into\NN$).
It is easy to see that a permutation of $\NN$ is completely determined
up to conjugacy by its \emph{cycle type}, that is, the number of
$n$-cycles for $n=1,2,\ldots,\infty$.  Since the cycle type can be
explicitly calculated and encoded as an element of the standard Borel
space $\NN^\NN$, we have that $\mathcal C_\NN$ is smooth.

If $E$ is not smooth, it is still very useful to speak of its
complexity relative to other classification problems.  If $E,F$ are
equivalence relations on standard Borel spaces $X,Y$, then we say that
$E$ is \emph{Borel reducible} to $F$, and write $E\leq_BF$, if there
exists a Borel map $f\from X\into Y$ such that for all $x,x'\in X$,
\[x\mathrel{E}x'\iff f(x)\mathrel{F}f(x')\;.
\]
If $E\leq_BF$, then the classification problem associated to $E$
should be regarded as no more difficult than that associated to $F$,
in the sense that any set of complete invariants for $F$ can be used
for $E$ as well.  Now, $\leq_B$ defines a partial preordering on all
equivalence relations on standard Borel spaces, and a foundational
goal of the subject is to understand the structure of this ordering.

It turns out that there is a $\leq_B$-maximum element among the
classification problems $\oiso_\sigma$ introduced earlier.  An
equivalence relation $E$ is said to be \emph{Borel complete} if for
every countable relational language $\mathcal L$ and every sentence
$\sigma$ of $\mathcal L_{\omega_1,\omega}$, we have
$\oiso_\sigma\leq_BE$.  A classification problem which corresponds to
a Borel complete equivalence relation should be regarded as totally
intractable.  The class of all countable structures is of course Borel
complete, but there are many more tangible Borel complete classes.
For instance, it is not difficult to see that the isomorphism relation
for countable graphs is Borel complete, and it is shown in \cite{FS}
that the classes of countable trees and linear orders are Borel
complete as well.

Now, we may in fact consider a class of equivalence relations which is
much broader than just the $\oiso_\sigma$.  An equivalence relation
$E$ is said to be \emph{classifiable by countable structures} if $E$
is Borel reducible to $\oiso_\sigma$ for some
$\mathcal{L}_{\omega_1,\omega}$-sentence $\sigma$ in a countable
relational language $\mathcal L$.  Thus, a Borel complete equivalence
relation is also universal for the class of equivalence relations
which are classifiable by countable structures.

A large source of examples of equivalence relations which are
classifiable by countable structures are those which arise from
certain Polish group actions.  Generally, if the Polish group $G$ acts
in a Borel fashion on the standard Borel space $X$, then the
corresponding orbit equivalence relation $E_G^X$ is an analytic
equivalence relation.  Here, $E_G^X$ is defined by
\[x\mathrel{E}_G^Xx' \iff \textrm{$x$ and $x'$ lie in the same $G$-orbit.}
\]
For instance, $\oiso_\sigma$ is the orbit equivalence relation induced
by the action of $S_\infty$ on $\Mod(\sigma)$.  By \cite{BK} and
\cite[Theorem~3.5.1]{gao}, any orbit equivalence relation induced by a
Borel action of a closed subgroup of $S_\infty$ is classifiable by
countable structures.

We shall now specialize to classification problems of the following
form.  If $\mathcal M$ is any countable model, then $\aut\mathcal M$
denotes its automorphism group, and $\mathcal{C}_{\mathcal{M}}$ the
conjugacy equivalence relation on $\aut\mathcal M$.  That is,
\[f\mathrel{\mathcal C}_{\mathcal M}g\iff(\exists h\in\aut\mathcal M)(hf=gh)\;.
\]
Then $f,g\in\aut\mathcal M$ are conjugate if and only if the
expansions $(\mathcal M,f)$ and $(\mathcal M,g)$ are isomorphic, and
so $\mathcal C_\Gamma$ is classifiable by countable structures.
(Alternatively, $\mathcal C_{\mathcal M}$ is the orbit equivalence
relation induced by the conjugation action of the closed subgroup
$\aut\mathcal M\leq S_\infty$ on itself, and hence is classifiable by
countable structures.)  We are particularly interested in the
complexity of $\mathcal C_{\mathcal M}$ for $\mathcal M$ a model of
some $\aleph_0$-categorical theory, as there are several interesting
examples in the recent literature.

\begin{thm*}[\protect{\cite[Theorem~76]{F}}]
  The conjugacy problem $\mathcal C_\QQ$ for the automorphism group of
  $(\QQ,\leq)$ is Borel complete.
\end{thm*}

\begin{thm*}[\protect{\cite[Theorem~5]{CG}}]
  Let $\mathbb B$ denote the countable atomless boolean algebra.  Then
  the conjugacy problem $\mathcal C_{\mathbb B}$ for the automorphism
  group of $\mathbb{B}$ is Borel complete.
\end{thm*}

We shall prove the analogous result for the automorphism group of the
random graph.  Recall that the \emph{random graph} $\Gamma$ is the
unique countably infinite graph satisfying the homogeneity property:
for any pair of finite and disjoint sets of vertices
$U,V\subset\Gamma$, there exists a vertex $x$ adjacent to each member
of $U$ and to no member of $V$.

\begin{thm}
  \label{thm_main}
  The conjugacy problem $\mathcal C_\Gamma$ for the automorphism group
  of the random graph $\Gamma$ is Borel complete.
\end{thm}

This result, which will be established in the next section, gives the
third of many conceivable examples of conjugacy problems for
automorphism groups of $\aleph_0$-categorical structures whose Borel
complexity is the maximum possible.  However, in many other cases the
conjugacy problem turns out to be smooth.  For instance, we have
already noted that the conjugacy problem for $S_\infty$ is smooth, and
it is not hard to see that the conjugacy problem for the automorphism
group of the complete binary tree is smooth as well.  This leaves open
the question of which Borel complexities can arise as $\mathcal
C_{\mathcal M}$ for some countable model $\mathcal M$ of an
$\aleph_0$-categorical theory.

\begin{conj*}
  If $\mathcal M$ is the countable model of an $\aleph_0$-categorical
  theory then $\mathcal C_{\mathcal M}$ is either smooth or Borel
  complete.
\end{conj*}

\section{Proof of the Main Theorem}

Our proof will proceed on lines similar to the proofs of Foreman's and
Camerlo-Gao's theorems.  The first step is to establish (or recall) the
Borel completeness of the isomorphism relation on the corresponding
class of structures.

\begin{thm*}[\protect{\cite[Theorem~3]{FS}}]
  The isomorphism relation on the space of countable linear orders is
  Borel complete.
\end{thm*}

\begin{thm*}[\protect{\cite[Theorem~1]{CG}}]
  The isomorphism relation on the space of countable Boolean algebras
  is Borel complete.
\end{thm*}

The second step is to reduce this isomorphism relation to the
conjugacy relation on the automorphism group of the universal
structure.  For instance, Foreman's argument proceeds as follows.
Given a countable linear ordering $x$ on $\NN$, construct an
automorphism $\phi_x$ of $\QQ$ whose fixed point set is isomorphic to
$x$, and whose nontrivial orbitals are all ``up-bumps.''  (Here, the
\emph{orbitals} of $\phi\in\aut\QQ$ are the convex closures of the
orbits $\set{\phi^n(q):n\in\ZZ}$; see \cite{F} for details).  This
construction can be carried out in such a way that if $x$ and $y$ are
isomorphic, then the linear orderings of orbitals of $\phi_x$ and
$\phi_y$ will be isomorphic.  A classical result (see for instance
\cite[Theorem~2.2.5]{G}) then implies that $\phi_x$ and $\phi_y$ are
conjugate. On the other hand, if $\phi_x$ and $\phi_{y}$ are
conjugate, then the linear orderings of fixed points of $\phi_x$ and
$\phi_y$ are clearly isomorphic, and hence $x$ is isomorphic to $y$.

In our case, we shall need the following:

\begin{thm*}[\protect{\cite[Theorem~1]{FS}}]
  The isomorphism relation on the space of countable, connected graphs
  is Borel complete.
\end{thm*}

Theorem~\ref{thm_main} then follows from the corresponding reduction.

\begin{thm}
  \label{thm_graphs}
  The isomorphism relation on the space of countable graphs is Borel
  reducible to $\mathcal C_\Gamma$.
\end{thm}

\begin{proof}[Proof of \protect{Theorem~\ref{thm_graphs}}]
  For this proof, we formally let $\Gamma$ denote a fixed copy of the
  random graph with underlying set $\NN$.  Here as usual we identify
  $\Gamma$ with its edge relation so that formally $\Gamma\in
  2^{\mathbb N\times\mathbb N}$.  Given a graph $x$ with underlying
  set $\NN$, we will build, in a Borel fashion, an automorphism
  $\phi_x$ of $\Gamma$ in such a way that:
  \begin{center}
    $x$ is isomorphic to $y$ $\iff$ $\phi_x$ is
    conjugate to $\phi_y$.
  \end{center}
  In order to do this, we shall first construct a graph
  $\Delta_x$ with the following properties:
  \begin{enumerate}
  \item $\Delta_x$ is isomorphic to the random graph $\Gamma$;
  \item $\Delta_x$ contains two isomorphic copies of $x$ as induced
    subgraphs; and
  \item there exists an automorphism of $\Delta_x$ which interchanges
    the two copies of $x$.
  \end{enumerate}
  This will be done using the ideas of \cite[Lemma~2.1]{M} and the
  comments that follow it.  Then, roughly speaking, we shall take
  $\phi_x$ to be the automorphism given by (c), thought of as an
  automorphism of $\Gamma$.

  Take for the vertices of $\Delta_x$ the set $\mathbb N\times\mathbb
  N$.  Visually, we think of the subset $\{i\}\times\mathbb N$ as the
  $i^\mathrm{th}$ \emph{row} of $\Delta_x$.  We shall put a copy of
  $x$ in each of rows $0,1$, and leave the remaining rows blank.  More
  precisely, let
  \[(i,m)\sim_{\Delta_x}(i,n)\iff m\sim_{x} n\qquad(\textrm{for }i=0,1)
  \]
  and if $i>1$, then $(i,m)\nsim_{\Delta_x}(i,m)$ for all
  $m,n\in\mathbb{N}$.

  We then define edges going across the rows as follows.  First, for
  each $n\in\mathbb N$ place an edge joining $(0,n)$ with $(1,n)$.
  For each row $i>1$, let $\seq{S^i_n:n\in\NN}$ be a (previously
  fixed) enumeration of the sets of $i^2$ vertices consisting of
  exactly $i$ vertices from each of rows $0,\ldots, i-1$.  Then place
  an edge from each vertex $(i,n)$ to each of the vertices in $S^i_n$.
  Thus, for each choice of $i$ vertices from each of rows
  $0,\ldots,i-1$, there is exactly one vertex in row $i$ which is
  adjacent to these vertices and to no others in rows $0,\ldots,i$.

  Note that the graph $\Delta_x$ satisfies the homogeneity property of
  the countable random graph, and hence is isomorphic to $\Gamma$.  In
  fact, we can build an explicit isomorphism by using a back and forth
  construction that chooses, at each stage in the construction, the
  least possible witness.  In this way we obtain a continuous
  assignment $x\mapsto\psi_x$, where for each graph $x\in
  2^{\NN\times\NN}$, $\psi_x\from\NN\into\NN$ is an isomorphism of
  $\Delta_x$ with $\Gamma$.

  Now, for any $x$, the map which interchanges the elements $(0,n)$
  and $(1,n)$ is a partial automorphism of $\Delta_x$.  Moreover, this
  map extends uniquely to a bijection $\phi\from\mathbb
  N^2\rightarrow\mathbb N^2$ which preserves every row $i\geq2$ of
  $\Delta_x$ setwise, and all of the edges added in the construction
  of $\Delta_x$.  Since $\phi$ is in particular an automorphism of
  $\Delta_x$, we may use the isomorphism $\psi_x$ between $\Delta_x$
  and $\Gamma$ to regard $\phi$ as an element of $\aut\Gamma$.
  Specifically, we let $\phi_x=\psi_x\,\phi\,\psi_x^{-1}$.  Since the
  map $x\mapsto\psi_x$ is Borel, the map $x\mapsto\phi_x$ is Borel
  too.

  We must now show that $x\cong y$ if and only if $\phi_x$ and
  $\phi_y$ are conjugate elements of $\aut\Gamma$.  Again following
  the isomorphisms $\psi_x$ and $\psi_y$, it suffices to show that
  $x\cong y$ if and only if there exists a graph isomorphism
  $\alpha\from\Delta_x\cong\Delta_y$ satisfying
  $\alpha\,\phi=\phi\,\alpha$.

  For the forward direction, if $\alpha_0\from x\cong y$ then
  $\alpha_0$ induces a row-preserving isomorphism
  $\alpha\from\Delta_x\cong\Delta_y$ which acts by $\alpha_0$ on row
  $0$ and row $1$ of $\Delta_x$.  It is easy to check that this
  $\alpha$ satisfies our requirements.

  Conversely, let $\alpha\from\Delta_x\rightarrow\Delta_y$ be a graph
  isomorphism satisfying $\alpha\,\phi=\phi\,\alpha$.  We wish to show
  that $\alpha$ somehow induces an isomorphism from $x$ to $y$.  We
  first show that $\alpha$ sends rows $0,1$ of $\Delta_x$ onto rows
  $0,1$ of $\Delta_y$.  Indeed, let $e\in\{0,1\}$ and $n\in\NN$, and
  suppose that $\alpha(e,n)=(i,\ell)$ where $i\geq2$.  Then we would
  have:
  \begin{align*}
  \alpha(1-e,n) &= \alpha\,\phi(e,n)\\
                &= \phi\,\alpha(e,n)\\
                &= \phi(i,\ell)\;.
  \end{align*}
  Since $\phi$ preserves row $i$ of $\Delta_y$, we conclude that
  $\alpha(1-e,n)=(i,\ell')$ for some $\ell'\neq \ell$.  Now, since
  $(e,n)\sim_{\Delta_x}(1-e,n)$, we must have
  $(i,\ell)\sim_{\Delta_y}(i,\ell')$, but this is a contradiction,
  since no two elements in row $i$ of $\Delta_y$ are adjacent.
  Similarly, $\alpha^{-1}$ maps rows $0,1$ of $\Delta_y$ to rows $0,1$
  of $\Delta_x$ and it follows that $\alpha$ maps rows $0,1$ of
  $\Delta_x$ bijectively onto rows $0,1$ of $\Delta_y$.

  We next show that $\alpha$ sends adjacent pairs $(0,n),(1,n)$ of
  $\Delta_x$ to adjacent pairs $(0,\ell),(1,\ell)$ of $\Delta_y$.
  Indeed, if $e,f\in\{0,1\}$ and $\alpha(e,n)=(f,\ell)$, then:
  \begin{align*}
    \alpha(1-e,n) &= \alpha\,\phi(e,n)\\
                  &= \phi\,\alpha(e,n)\\
                  &= \phi(f,\ell)\\
                  &= (1-f,\ell)\;.
  \end{align*}
  Note that we are not claiming that $\alpha$ need send the whole of
  row $0$ to the whole of row $0$ or $1$, and indeed this need not be
  the case.  But it follows from what we have shown that $\alpha$
  induces a bijection $\alpha_0\from x\rightarrow y$ defined by
  \[\alpha(e,n)=(f,\alpha_0(n))
  \]
  for some $e,f\in\{0,1\}$.

  It remains only to show that $\alpha_0$ is a graph isomorphism.
  Indeed, if $m\sim_x n$, then $(0,m)\sim_{\Delta_x}(0,n)$ and so
  $\alpha(0,m)\sim_{\Delta_y}\alpha(0,n)$.  Hence, there are
  $e,f\in\set{0,1}$ such that
  $(e,\alpha_0(m))\sim_{\Delta_y}(f,\alpha_0(n))$.  Since
  $\alpha_0(m)\neq \alpha_0(n)$, by the construction of $\Delta_y$ we
  must have that $e=f$, and hence that
  $\alpha_0(m)\sim_{y}\alpha_0(n)$.  This shows that $\alpha_0$ is a
  graph embedding.  Moreover, we can repeat this argument with
  $\alpha_0^{-1}$ to conclude that $\alpha_0$ is a graph isomorphism.
\end{proof}

\bibliographystyle{alpha}
\begin{singlespace}
  \bibliography{summer}
\end{singlespace}

\end{document}